%% file: mow2016.tex
\documentclass{article}

\usepackage{multirow,lscape}

\usepackage{amsmath}
\usepackage{amsfonts}
\usepackage{amssymb}
\usepackage{amsthm}
\usepackage{todonotes}
\usepackage{subfiles}
\usepackage{soul}
\usepackage{multicol}
\usepackage{tabularx}

\newtheorem{proposition}{Proposition}
\newtheorem{assumption}{Assumption}

\newcommand\Sos{{\operatorname{\text{\small $\Sigma$}}}}
\newcommand{\abs}[1]{\left|#1\right|}

\usepackage{color}

\begin{document}
\title{Adjustably Robust Optimal Power Flows with Demand Uncertainty via Path-based Flows}

\author{
Jakub Marecek$^1$, Adam Ouorou$^2$, Guanglei Wang$^2$\thanks{Guanglei Wang is the first author, despite being listed the last
in the alphabetical ordering. Email: guanglei.wang@orange.com}\\[3mm]
$^1$ {\small IBM Reserach -- Ireland, Dublin, Ireland}\\
$^2$ {\small Orange Labs Research, Issy-les-Moulineaux, France}
}

\maketitle

\begin{abstract}
We study the optimal power flow problem with switching (or, equivalently, the line expansion problem) under demand uncertainty.
Specifically, we consider the line-use variables at the first stage and the current- or power-flow at the second stage of two affinely adjustably robust formulations. 
\end{abstract}

\subfile{dcopf-aar-introduction}

\subfile{dcopf-aar-notation}
\newpage
\subfile{dcopf-aar-deterministic}
\subfile{dcopf-aar-robust}
\subfile{dcopf-aar-pop}
\subfile{dcopf-aar-approx}

\begin{figure}[t]
\centering
\includegraphics[width=0.99\textwidth,clip=true,trim=4.5cm 20cm 4cm 4.5cm]{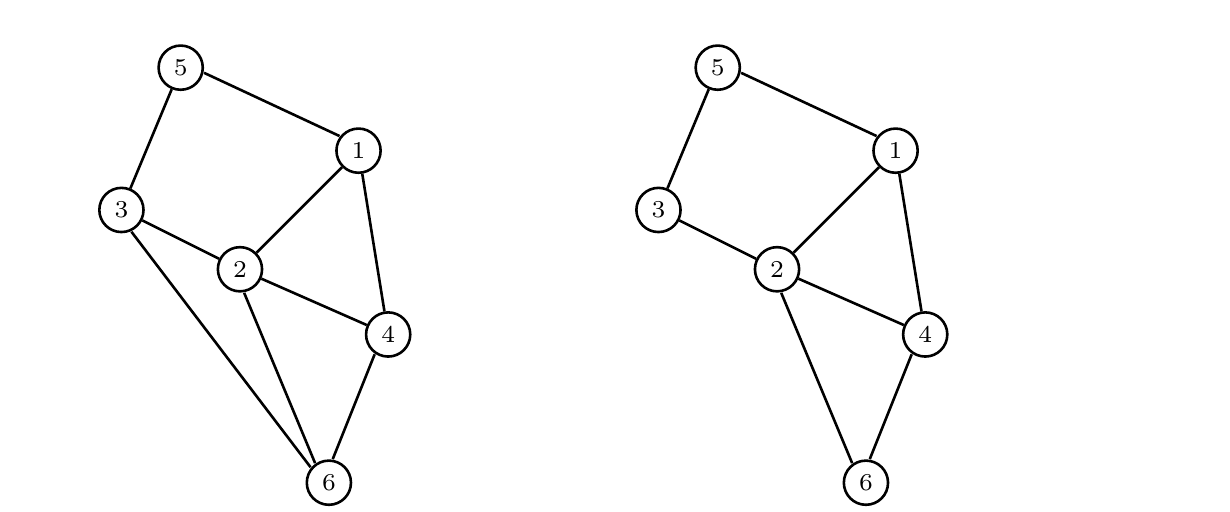}
\caption{The instance Garver6y with the globally optimal solution (left) and 
the only other AC-feasible solution (right).
}
\label{fig:garver6instance}
\end{figure}

\subfile{dcopf-aar-results}

\section{Conclusions}
We have taken first steps toward the study of adjustable robustness 
with respect to demand uncertainty in power systems. 
Much work remains to be done, starting with extensive computational testing,
both in terms of the comparison the solutions to other methods for dealing with uncertainty, and in terms of scalability of run-time.
Further, one should like to study the impact of losses in alternating-current
models on the solutions.
First illustrations in another paper submitted to this conference \cite{JonasPSCC}, which considers two-stage stochastic programs with mixed-integer polynomial optimisation recourse, suggest the impact is considerable.
This seems to open up a large area for future research.


\begin{landscape}
\begin{table*}[th!]
\centering
Buses:\\[2mm]
\begin{tabular}{|l|l|r|r|r|r|r|r|r|r|r|r|r|}\hline
bus & type &   $P_d$   &   $Q_d$  &  $G_s$  &  $B_s$  & area &  $V_m$  &  $V_a$  & baseKV & zone & $V_{\max}$ & $V_{\min}$ \\ \hline 
1 & 3 &  80.00 & 16.00 & 0.00 & 0.00 & 1.00 & 1.00 & 0.00 & 230.00 & 1.00 & 1.05 & 0.95 \\ \hline 
2 & 1 & 240.00 & 48.00 & 0.00 & 0.00 & 1.00 & 0.00 & 0.00 & 230.00 & 1.00 & 1.05 & 0.95 \\ \hline 
3 & 2 &  40.00 &  8.00 & 0.00 & 0.00 & 1.00 & 1.00 & 0.00 & 230.00 & 1.00 & 1.05 & 0.95 \\ \hline 
4 & 1 & 160.00 & 32.00 & 0.00 & 0.00 & 1.00 & 0.00 & 0.00 & 230.00 & 1.00 & 1.05 & 0.95 \\ \hline 
5 & 1 & 240.00 & 48.00 & 0.00 & 0.00 & 1.00 & 1.00 & 0.00 & 230.00 & 1.00 & 1.05 & 0.95 \\ \hline 
6 & 2 &  0.00  &  0.00 & 0.00 & 0.00 & 1.00 & 0.00 & 0.00 & 230.00 & 1.00 & 1.05 & 0.95 \\ \hline 
\end{tabular}\\[6mm]
Generators:\\[2mm]
\begin{tabular}{|c|c|c|c|c|c|c|c|c|c|c|c|c|c|}\hline
bus &  Qmax  &  Qmin  &  Vg  &  Pmax  & Pmin &  Pc1 &  Pc2 & Qc1min & Qc1max & Qc2min & Qc2max \\ \hline 
1 & 48.25 & -10.00 & 1.00 & 160.00 & 0.00 & 0.00 & 0.00 &  0.00  &  0.00  &  0.00  &  0.00  \\ \hline 
3 & 101.25 & -10.00 & 1.00 & 370.00 & 0.00 & 0.00 & 0.00 &  0.00  &  0.00  &  0.00  &  0.00  \\ \hline 
6 & 183.00 & -10.00 & 1.00 & 610.00 & 0.00 & 0.00 & 0.00 &  0.00  &  0.00  &  0.00  &  0.00  \\ \hline 
\end{tabular}\\[6mm]
Branches:\\[2mm]
\begin{tabular}{|l|l|r|r|r|r|r|r|r|}\hline
fbus & tbus & r & x & rateA  & rateB & rateC &	ratio & angle \\  \hline 
1 & 2 &   0.040 &  0.40 &   180 &  250 & 250 & 0 & 0 \\ \hline 
1 & 4 &   0.060 &  0.60 &   150 &  250 & 250 & 0 & 0 \\ \hline 
1 & 5 &   0.010 &  0.10 &   360 &  250 & 250 & 0 & 0 \\ \hline 
2 & 3 &   0.020 &  0.20 &   180 &  250 & 250 & 0 & 0 \\ \hline 
2 & 4 &   0.040 &  0.40 &   180 &  250 & 250 & 0 & 0 \\ \hline 
2 & 6 &   0.015 &  0.15 &   360 &  250 & 250 & 0 & 0 \\ \hline 
3 & 5 &   0.010 &  0.10 &   360 &  250 & 250 & 0 & 0 \\ \hline 
3 & 6 &   0.024 &  0.24 &   360 &  250 & 250 & 0 & 0 \\ \hline 
4 & 6 &   0.008 &  0.08 &   360 &  250 & 250 & 0 & 0 \\ \hline
\end{tabular}
\caption{The details of the instance depicted in Figure~\ref{fig:garver6instance} in the Matpower format. Columns not listed are uniformly at the default values.}
\label{tab:1}
\end{table*}
\end{landscape}

\bibliographystyle{abbrv} 
\bibliography{acopf,aar,literature}
\clearpage
\subfile{dcopf-aar-appendix}

\end{document}

%% file: dcopf-aar-introduction.tex
\section{Introduction}

There is a long history of research into optimisation under uncertainty  in the
operations of power systems. Often, one considers a small discrete number of
realisations of uncertainty within the framework of stochastic programming.
At other times, one should like to consider an infinite number of realisations of
uncertainty, e.g., by replacing uncertain scalar values with uncertainty sets, such as intervals. In robust optimisation, each \textit{uncertainty-immunized solution} is required to be feasible for all the realisations of the uncertainty. The robust counterpart of
an uncertain problem seeks such a solution and can often be solved efficiently. This approach has a two-fold appeal: first, the robust counter-parts are often easier 
to solve than a scenario-expansion in stochastic programming,
and second, one can often derive guarantees as to the performance of the solutions.

The applications of robustness in power systems have, however, been limited due to the difficulties in dealing with equalities in robust optimisation. 
Even very recent work considers only 
inequality constraints, rather than equalities required in transmission-constrained problems. 
One could use two-level formulations \cite{Bentaletal2004},
also known as adjustable robustness or the decision rules approach,
but has to tackle the resulting non-convex non-linear optimisation problems.

We propose dealing with expansion decisions (line-use variables) 
$x_l$ for branches $l\in E$ at the first, non-adjustable stage. The adjustable variables are power flow variables $s$ and
power reservation variables $p$. 
Demand 
$\xi_k$ at bus $k\in \mathcal{K}$ is unknown a priori, but 
estimates in one of two forms are available.
The first form is based on the nominal value, the dispersion range
around the nominal value, and the number of nomila values, which are
allowed to deviate from the nominal value, in the spirit of 
{\it budget uncertainty} \cite{Bentaletal2004}.
The second form could be seen either as based on the value at risk, 
or as a polyhedral approximation of certain ellipsoidal uncertainty sets.


In order to solve the problem, we introduce a number of novel techniques:
\begin{itemize}
\item a power-flow formulation based on path-based variables for current or power
and Kirchhoff's first law (Section~\ref{sec:deterministic})
\item an elaborate construction of adjustably robust counter parts 
considering little studied uncertainty sets (Section ~\ref{sec:robust})
\item methods from polynomial optimization \cite{Ghaddar2015} for solving
the robust counter parts exactly (Section~\ref{sec:computational})
\item for comparison, an approximation of adjustable decision rules by affine decision rules, which amounts to requiring $s$ and $p$ to be of the form 
\begin{equation} \label{ADR}
    \begin{split}
    p_g(\xi) &= p_g^0 + \sum\limits_{h\in\mathcal{K}}p_g^h\xi_h, \hspace{0.5em}  g\in G \\
    s_p(\xi) &= s_p^0 + \sum\limits_{h\in\mathcal{K}}s_p^h\xi_h, \hspace{0.5em}  p\in \mathcal{P}_k, k\in \mathcal{K} \\
    \end{split}
\end{equation}
for buses $\mathcal{K}$, paths $\mathcal{P}_k$ to bus $k\in \mathcal{K}$, and 
generators $G$ (Section~\ref{sec:approx}). 
\end{itemize}
A discussion of the options is presented in conclusions.
This is complementary to recent work on 
adjustable robustness \cite{Jabr2013OPF,Jabr2014MOPF}, 
robustness in terms of conditional value at risk \cite{Zhang2013} 
and chance-constrained models \cite{Bienstock2013}.
It may suggest an appealing direction for further research, e.g. making it possible to derive guarantees as to the performance of the solutions \cite{bandi2012}.

%% file: dcopf-aar-notation.tex
\subsection{Notation}

 \begin{tabularx}{0.99\textwidth}{lX}
     \multicolumn{2}{l}{Sets} \\
     $N$                                         & set of buses in the power network. \\
     $L$                                         & set of transmission lines.\\
     $G$                                         & set of generators $G\subseteq N$.\\
     $\mathcal{K}$                               & set of demands,    $\mathcal{K}\subseteq N$.\\
     $\mathcal{P}_k $                            & set of admissible paths from any generator \\
                                                 & to the customer $k \in \mathcal{K}$.\\
     $\mathcal{P}_g $                            & set of admissible paths from generator $g \in \mathcal{G}$ to any customer\\
     $\mathcal{P}_{kg}$                          & set of admissible paths between demand $k$ 
                                                 and generator $g$, where $ (k,g)\in \mathcal{K} \times \mathcal{G}$.\\
      \multicolumn{2}{l}{Parameters} \\
     $ f_l            $                          & the investment costs for line   $ l\in L  $\\
     $ R_l            $                          & the resistance of line $ l \in L $\\
     $ B_l            $                          & the susceptance of line $ l \in L $\\
     $ c_g            $                          & per-unit costs for generator $ g\in G  $\\
     $ d_k          $                            & the demand of customer $ k\in \mathcal{K} $\\
     $ \overline{y_l} $                          & an upper bound for the current flow for $l \in L$\\
     \multicolumn{2}{l}{Variables} \\
     $p_g$                                       & reserve power capacity for generator $g\in G$\\
     $s_l$                                       & current flow along line $ l \in L$ \\
     $x_l$                                       & indicator of switching line $ l \in L$ (0 open)\\
     \multicolumn{2}{l}{Uncertainty sets} \\
     $\Xi$                                         & budget-of-uncertainty u. set \eqref{defineXi} \\
$\bar d_k$ & nominal value of demand \\
$\hat d_k$ & dispersion value of demand \\
$\xi_k$ & a random variable in $[-1, 1]$ \\
$\kappa$ & parameter regulating dispersion \\
$\tau$ & number of deviations from nominal demands\\
    $\Xi'$                                       & value-at-risk uncertainty set \eqref{defineXiPrime} \\
    $\alpha$                                       & quantile for the ``$\alpha$-tail''
\end{tabularx}

%% file: dcopf-aar-deterministic.tex
\section{A Path-Based Formulation}
\label{sec:deterministic}
We propose a formulation, where for each customer, there is a scalar decision
variable for each path from a generator to the particular customer, and for
each edge on that path, a thermal limit constraint across all paths using it.

The objective is to minimize the transmission expansion cost $\sum_{l \in E
}f_lx_l$ and power reservation cost $\sum_{g\in G}  c_gp_g$. In this
paper, transmission expansion assumes, that the existing network is fixed
and we can only expand the network capacity by adding parallel transmission
lines between existing towers. 

Based this path-based formulation, we propose two variants. Formally, variable $y_{p}$ captures the current flow along path $p$.
\begin{align}
     \min& \hspace{0.5em} \sum\limits_{l \in E }f_lx_l + \sum\limits_{g\in G}  c_gp_g \tag{PB1}\label{pb1} \\
    \mathrm{s.t.}& \hspace{0.5em} \sum\limits_{p\in \mathcal{P}_k} y_{p} \geq
    \sqrt{ d_k / R }, & k&\in \mathcal{K},\\
                      &\hspace{0.5em} \sum\limits_{p \in \mathcal{P}_g}  y_p  \leq
    \sqrt{p_g/R}, &g& \in G , \\
    &\hspace{0.5em}  \sum\limits_{k\in \mathcal{K}}\sum\limits_{p\in\mathcal{P}_k: l\in p}
    y_{p}  \le \bar{y_l}x_l, &l& \in E, \\
                                  &\hspace{0.5em}  y_{p} \geq 0, p_g \geq 0, x_l \in \{0, 1\} ,
    &k&\in \mathcal{K}, p\in \mathcal{P}_k, g\in G, l\in E. 
\end{align}
Alternatively, variable $s_{p}$ captures the apparent power provided to customer $k$ along path $p$:
\begin{align}
    \min               & \hspace{0.5em}  \sum\limits_{l \in E }f_{l}x_{l} + \sum\limits_{g\in G}    c_g p_g     \tag{PB2} \label{pb2}\\
    \mathrm{s.t.}     & \hspace{0.5em} \sum\limits_{p\in \mathcal{P}_k} s_{p} \ge  d_k, & k & \in \mathcal{K},\label{demand}\\
                                          & \hspace{0.5em} \sum\limits_{p \in \mathcal{P}_g} s_{p}  \leq p_g, & g & \in G, \label{generation}\\
                                          & \hspace{0.5em} \sum\limits_{k\in \mathcal{K}}\sum\limits_{p\in\mathcal{P}_k: l\in p}
    s_{p} \le  \bar{y_l}^2R_lx_l,    & l& \in E,\label{thermal} \\
                                     & s_{p} \geq 0, p_g \geq 0, x_l\in \{0,1\} ,& k & \in \mathcal{K}, p\in\mathcal{P}_k, g\in G, l\in E.\label{domain}
\end{align}
Our work on this path-based formulation was partially motivated  by the work
of \cite{Kocuk2014}, who have shown, that one may formulate the optimal power
flow in the directed current model using Kirchoff's second law (Kirchoff's
Voltage Law), which requires the voltage differences along any {\it directed}
cycle $C$ to sum to zero. By the fact that any directed cycle $C: s-t-s$ can be
decomposed into two paths $P_{st}$ and $P_{ts}$, the following proposition holds.  
\begin{proposition}
The objective function values of path-based OPF and cycle-based OPF are equivalent. 
\end{proposition}

%% file: dcopf-aar-robust.tex
\section{The Adjustable Robust 
Formulation under Demand Uncertainty}
\label{sec:robust}


In the remaining of this paper, we derive the adjustable robust counterparts for model
\eqref{pb2}. In this paper, the uncertainty comes from the prediction of customer demands
or loads. Instead of trying either to approximate one 
particular multi-variate distribution of the load, as in stochastic programming,
or trying to work with any distributions on the particular support, as in early
work on robust optimisation, we consider two non-trivial uncertainty sets. 

In order to construct the first uncertainty set $\Xi$, we express the prediction 
of demands as $d_k = \bar d_k +
\xi_k\hat d_k$, where $\bar d_k, \hat d_k$ are nominal value and dispersion
value and $\xi_k$ is a random variable in $[-1, 1]$.  Following the notion of
{\it budget uncertainty} \cite{Bertsimas2004,Bertsimas2011,Bentaletal2004},  we define $\Xi$ as: 
\begin{align}
\label{defineXi}
\Xi=\left\{\xi\in\mathbb{R}^{\left|\mathcal{K}\right|}:\xi_i\in [-1,1], i\in
\mathcal{K},\left \|\xi\right\|_1 \le\kappa,\left\|\xi\right\|_\infty
\le\tau\right\}
\end{align}
where the values of $\kappa, \tau$ are usually determined by the decision maker
who is responsible for the outcomes of the set $\Xi$. The larger the values, the larger
the set, and the higher robust price. If $\kappa= 1, \tau = \left|\mathcal{K}
\right|$, then all the possible outcomes are considered; while if both values
are zero, the values of demands are just nominal values. In this case, the
problem can be reduced to the deterministic models \eqref{pb1} or \eqref{pb2}. 
Further, we consider an alternative uncertainty set, 
\begin{align}
\label{defineXiPrime}
    \Xi'=\left\{\boldsymbol{d} \in\mathbb{R}^{\left|\mathcal{K}\right|}: \right. & \exists \lambda \in
        \mathbb{R}^N_+,  \boldsymbol{1}\lambda = 1, \\ 
& \left. \boldsymbol{d} =
\sum\limits_{i=1}^N \lambda_i d_i, \lambda_i \le \frac{1}{N(1-\alpha)} \right\}, \notag
\end{align}
where $\alpha$ is determined by the decision maker. It assumes that $N$
observations are available and controls the ``$\alpha$-tail'' of the uncertainty.
This seems to have been considered only once before, see page 467 of \cite{Bertsimas2011}.
The set $\Xi'$ can be seen as a polyhedral approximation of the ellipsoidal set considered
by Ben-tal et al.\cite{Bentaletal2004}.  

Subsequently, we can construct the adjustable robust counterpart of the
uncertain program for \eqref{pb2} with $\Xi$ and $\Xi'$. In both cases,  the
non-adjustable decision variables are expansion variables
$x_l, l\in E$. The adjustable variables are power flow variables $s$ and
power reservation variables $p$. Hence, we formulate the problem as 
\begin{align}
    \min          &\hspace{0.5em} \sum\limits_{l \in E }f_{l}x_{l}
    +\max\limits_{\xi\in \Xi}\min\limits_{s(\xi),p(\xi)}  \sum\limits_{g\in G}    c_g p_g  \tag{PB2-ARC} \label{pb2-arc}\\
    \mathrm{s.t.} & \hspace{0.5em} x \in \bigcap\limits_{\xi\in \Xi}\boldsymbol{Proj}_x \mathcal{F}(\xi) 
\end{align}
where, 
\begin{align*}
    \mathcal{F}(\xi) =\left \{ 
        (x,p(\xi),s(\xi)): \textrm{ parameterized }\eqref{demand},  \eqref{generation}, \eqref{thermal}, \eqref{domain} ~\textrm{hold}
    \right\}.
\end{align*}
We then write out the dual of the inner minimization problem with dual
associated multipliers $\lambda,\phi,\eta$ and denote by $S(x,\xi)$ the power
generation cost.
{\small
\begin{align*}
 S(x,\xi) = \max \limits_{\lambda, \phi,\eta}&\hspace{0.5em}\sum\limits_{k\in
\mathcal{K}} \lambda_kd_k - \sum\limits_{l\in E} \bar{y_l}^2x_l R_l \eta_l \\
               \textrm{s.t.}&\hspace{0.5em}  c_g -  \phi_g \geq 0,   \forall g\in G\\
                            &\hspace{0.5em}\sum\limits_{l\in E: l\in p}\eta_l - \lambda_k + \phi_g \geq 0,
            \forall (k,g)\in \mathcal{K}\times G, \forall p\in\mathcal{P}_{kg}\\
                         &\hspace{0.5em} \lambda, \phi,\eta \geq 0.
\end{align*}
}
By strong duality, the inner minimization problem is equivalent to
$\max\limits_{\xi\in
\mathcal{U}} S(x,\xi) $. Furthermore,  \eqref{pb2-arc} might be reduced to
 the following bilinear program:
 {\small
\begin{align}
\min  &\hspace{0.5em} \gamma +  \sum\limits_{l \in E }f_{l}x_{l} \tag{PB2-PP}\label{PB2-PP}\\
    \textrm{s.t.}&\hspace{0.5em}  c_g - \phi_g \geq 0,  g\in G, \label{PB2.1}\\
 &\hspace{0.5em}\sum\limits_{l\in E: l\in p}\eta_l - \lambda_k + \phi_g \geq
    0,\forall (k,g)\in \mathcal{K}\times G,p\in\mathcal{P}_{kg},\label{PB2.2}\\
 &\hspace{0.5em} \sum\limits_{k\in \mathcal{K}} (\lambda_k \bar d_k +\hat
    d_k\lambda_k\xi_k) - \sum\limits_{l\in  E} \bar{y_l}^2R_lx_l\eta_l \leq
    \gamma, \xi \in \Xi, \label{collection}\\
   &\hspace{0.5em} \lambda_k, \phi_g,\eta_l \geq 0, x_l\in \{0, 1\}, \forall
    k,g,l.\label{PB2.4}
\end{align}
}
To deal with a infinite number of of constraints \eqref{collection}, the robust
counterpart of constraints \eqref{collection} can be expressed as follows,
\begin{align}
\sum\limits_{k\in \mathcal{K}} \lambda_k\bar d_k - \sum\limits_{l\in  E}
\bar{y_l}^2R_lx_l\eta_l + \kappa t + \sum\limits_{k\in \mathcal{K}} \tau w_k
\le \gamma, \label{PB2xi1}\\
w_k+t+ \lambda_k\hat d_k \ge 0 , k\in \mathcal{K},\label{PB2xi2}\\
w_k+t- \lambda_k\hat d_k  \ge 0, k\in \mathcal{K},\label{PB2xi3}\\
w_k \ge 0, t\ge 0,  k\in \mathcal{K}\label{PB2xi4}.
\end{align}
Alternatively, for the uncertainty set $\Xi'$ 
we can express the robust counterpart of constraints \eqref{collection} as,
\begin{align}
    t+\frac{1}{N(1-\alpha)} \sum\limits_{i=1}^N w_i -\sum\limits_{l\in
E} \bar{y_l}^2R_lx_l\eta_l  \le \gamma,\label{PB2x1} \\
t + w_i \geq \sum\limits_{k\in \mathcal{K}} \lambda_kd_k^i  , i=1, \dots, N,\label{PB2x2}\\
    w_i \ge 0 ,i=1, \dots, N.\label{PB2x3}
\end{align}
To conclude this section, we present the exact adjustable robust counterpart
with uncertainty set $\Xi$ 
\begin{align}
    \min  &\hspace{0.5em} \gamma +  \sum\limits_{l \in E }f_{l}x_{l}
    \tag*{[PB2-ARC$\Xi$]}\\
    \textrm{s.t.}&\hspace{0.5em} \eqref{PB2.1},
    \eqref{PB2.2},\eqref{PB2.4},\eqref{PB2xi1},\eqref{PB2xi2},\eqref{PB2xi3},\eqref{PB2xi4}
\end{align}
and  the exact adjustable robust counterpart with uncertainty set $\Xi'$:
\begin{align}
    \min  &\hspace{0.5em} \gamma +  \sum\limits_{l \in E }f_{l}x_{l}  \tag*{[PB2-ARC$\Xi'$]}\\
    \textrm{s.t.}&\hspace{0.5em} \eqref{PB2.1},
    \eqref{PB2.2},\eqref{PB2.4},\eqref{PB2x1},\eqref{PB2x2},\eqref{PB2x3}
\end{align}

%% file: dcopf-aar-pop.tex
\section{Computations}
\label{sec:computational} 
The two-stage robust counterpart is NP-Hard to solve, in general 
\cite{Bentaletal2004,Bertsimas2011}.
Nevertheless, 
we derive one hierarchy of
semidefinite programming relaxation due to Lasserre \cite{Lasserre2006,Handbook},
which are asymptotically convergent, in theory, and solvable, in practice.

Let us sketch the approach out, following the treatment of \cite{Ghaddar2015}, on a generic polynomial optimisation problem
\begin{align}
\min \quad& f(x) \notag \\
\mbox{s.t.  }& g_i(x) \geq 0 \qquad i=\{1,\dots,m\} \tag*{[PP]}
\end{align}
If there exist homogeneous polynomials of degree $d$ in $n$-dimensional vector $x$, denoted $g_1(x),\ldots,g_k(x)$ such that 
$h(x)=\sum_{i=1}^k g_i(x)^2$,
polynomial $h(x)$ of degree $2d$ is sum-of-squares (SOS).
We use 
$\mathcal{P}_d(S)$ to denote the cone of polynomials of degree at most $d$ that are non-negative over some $S\subseteq \mathbb{R}^n$ and
$\Sigma_d$ to denote the cone of polynomials of degree at most $d$ that are SOS.
Choi et al.\ \cite{Choi1995} showed that each $\mathcal{P}_{2d}(\mathbb{R}^n)$ can be approximated as closely as desired by a sum-of-squares of polynomials.
Lasserre reformulated [PP] as
\begin{align}
\max \quad & \varphi & \mbox{s.t.  }& f(x)-\varphi  \geq 0 \quad \forall \: x \in S_G, \notag \\
= \max \quad & \varphi & \mbox{s.t.  } &f(x)-\varphi \in \mathcal{P}_d(S_{\mathcal{G}}). \tag*{[PP-D]}
\end{align}
where $\mathcal{G}=\{g_i(x): i=1,\dots,m \}$ and $S_{\mathcal{G}}=\{x \in \mathbb{R}^n :  g(x) \geq 0, \; \forall g \in \mathcal{G}\}$ 
In \cite{WKKM,Lasserre2006,KojimaMuramatsu2009}, the $n\times n$ matrix:
$$\mathcal{R}_{ij} = \begin{cases} 
\star& \text{for } i=j\\ 
\star& \text{for } x_i, x_j \text{ in the same monomial of } f \\ 
\star& \text{for } x_i, x_j \text { in the same constraint } g_k \\ 
0 & \text {otherwise},
\label{eq:CSP}
\end{cases} $$
and its associated adjacency graph $G$
represent the so called correlative sparsity of [PP].
Let $\{I_k \}_{k=1}^p$ be the set of maximal cliques of a chordal extension of $G$ following the construction in \cite{WKKM}, i.e. $I_k \subset \{ 1,\ldots,n \}$. 
This leads to approximation $\mathcal{{K}}^{r}_{\mathcal{G}}(I)$ of $\mathcal{P}_d(S)$:
\begin{align*}
  \mathcal{{K}}^{r}_{\mathcal{G}}(I) = \sum_{k=1}^p \left(\Sigma_r(I_k) +  \sum_{j\in J_k} g_j \Sigma_{r-\deg(g_j)}(I_k)\right),
\end{align*}
where $\Sigma_d(I_k)$ is the set of all sum-of-squares polynomials of degree up to $d$ supported on $I_k$ and $(J_1,\ldots, J_p)$ is a partitioning of the set of polynomials $\{g_j\}_j$ defining $S$ such that for every $j$ in $J_k$, the corresponding $g_j$ is supported on $I_k$. The support $ I \subset \{1,\ldots,n\} $ of a polynomial contains the indices $i$ of terms $x_i$ which occur in one of the monomials of the polynomial.  The sparse hierarchy of SDP relaxations is then given by 
{\small\begin{align}
&\max_{\varphi,\sigma_k(x), \sigma_{r,k}(x)} \:   \varphi \tag*{[PP-SH$_r$]$^*$} \label{eq-Lass-sparse}  \\ 
&\mbox{s.t. } f(x)-\varphi=  \sum_{k=1}^p \left(\sigma_{k}(x) +  \sum_{j\in J_k} g_j(x) \sigma_{j,k}(x)\right) \notag \\
&\quad \sigma_{k}\in\Sos_r((I_k)), \sigma_{j, k} \in \Sos_{r-\deg(g_j)}(I_k). \notag 
\end{align}}
We denote the dual of [PP-SH$_r$]$^*$ by [PP-SH$_r$], in keeping with \cite{Ghaddar2015}.

Following Lasserre \cite{Lasserre2006}, we introduce:

\begin{assumption}
\label{as:Lasserre}
Let $S$ denote the feasible set of a problem of form [PP]. Let $\{I_k \}_{k}$ denote the $p$ maximal cliques of a chordal extension of the sparsity pattern graph of the [PP].
\begin{enumerate}
\item[(i) ]
Then, there is a $M>0$ such that $\parallel x \parallel_{\infty} < M$ for all $x\in S$.
\item[(ii) ]
Ordering conditions for index sets as in Assumption 3.2 (i) and (ii) in \cite{Lasserre2006}.
\item[(iii) ]
Running-intersection-property, c.f. \cite{Lasserre2006}, holds for $\{I_k\}_k$.
\end{enumerate}
\end{assumption}

Notice that this assumption is easy to satisfy for our robust counter parts, where 
all variables are bounded from above and below, and $S$ is hence compact. (One can add redundant quadratic inequalities such that Assumption \ref{as:Lasserre} (i) is satisfied; (ii) and (iii) can be satisfied by construction of the sets $\{I_k \}_{k}$ and the chordal extension, as pointed out in \cite{KojimaMuramatsu2009}.)
Indeed, the resulting full-dimensional polynomial optimisation problems are very nicely behaved.

Following Lasserre \cite{Lasserre2006} further, we have:

\begin{proposition}[Asymptotic Convergence]
\label{thm:convergent}
Whenever Assumption \ref{as:Lasserre} holds for the feasible sets of [PB2-ARC$\Xi$] and [PB2-ARC$\Xi'$], for the sparse hierarchy [PB2-ARC$\Xi$-SH$_r$]$^*$ for [PB2-ARC$\Xi$] and [PB2-ARC$\Xi'$-SH$_r$]$^*$ for [PB2-ARC$\Xi'$] and their respective duals holds: 
\begin{enumerate}
\item[(a) ]
$\inf \text{[PB2-ARC$\Xi$-SH$_r$]} \nearrow \min \text{([PB2-ARC$\Xi$])} \text{ as } r\rightarrow\infty,$\\
\item[(b) ]
$\sup \text{[PB2-ARC$\Xi$-SH$_r$]$^*$} \nearrow \min \text{([PB2-ARC$\Xi$])} \: \text{as} \: r\rightarrow\infty,$\\
\item[(c) ]
If the interior of the feasible set of [PB2-ARC$\Xi$] is nonempty, there is no duality gap between [PB2-ARC$\Xi$-SH$_r$] and [PB2-ARC$\Xi$-SH$_r$]$^*$.\\
\item[(d) ]
If the interior of the feasible set of [PB2-ARC$\Xi$] is nonempty, 
there is no duality gap between [PB2-ARC$\Xi$-SH$_1$] and [PB2-ARC$\Xi$-SH$_1$]$^*$.
\item[(e) ] 
If [PB2-ARC$\Xi$] has a unique global minimizer, then as $r$ tends to infinity the components of the optimal solution of [PB2-ARC$\Xi$-SH$_r$] corresponding to the linear terms converge to  
the unique global minimiser of [PB2-ARC$\Xi$].
\end{enumerate}
Throughout, the analogous result holds for [PB2-ARC$\Xi'$].
\end{proposition}
\begin{proof}
The proof is the same as in \cite{Ghaddar2015}: (a,b,e) follow from Theorem~3.6 of  \cite{Lasserre2006}, (c) follows from Theorem 5 of \cite{KojimaMuramatsu2009}, and (d) from  \cite{WKKM}.
\end{proof}

%% file: dcopf-aar-approx.tex
\section{Affine decision rule approximation} \label{sec:approx}
As a comparison with the exact formulations [PB2-ARC$\Xi$] and [PB2-ARC$\Xi'$],
we approximate the adjustable decision rules by the affine decision rule on the
demands, which amounts to requiring $s(\xi)$ and $p(\xi)$ to be of the form
\begin{equation} \label{ADR}
    \begin{split}
    p_g(\xi) &= p_g^0 + \sum\limits_{h\in\mathcal{K}}p_g^h\xi_h, \hspace{0.5em}  g\in G \\
    s_p(\xi) &= s_p^0 + \sum\limits_{h\in\mathcal{K}}s_p^h\xi_h, \hspace{0.5em}  p\in \mathcal{P}_k, k\in \mathcal{K} \\
    \end{split}
\end{equation}
for buses $\mathcal{K}$, paths $\mathcal{P}_k$ to bus $k\in \mathcal{K}$, and 
generators $G$. With these affine functions, we rewrite the affinely adjustable robust counterpart
\eqref{pb2-arc} as
\begin{align}
    \min          &\hspace{0.5em} \sum\limits_{l \in E }f_{l}x_{l} + \gamma  \tag{PB2-AAR} \label{pb2-aar}\\
    \mathrm{s.t.} &  \sum\limits_{g\in G}  c_g
    (p_g^0+\sum\limits_{h\in\mathcal{K}}p_g^h\xi_h) \leq \gamma, ~\xi\in\Xi,\\
                  &\sum\limits_{p\in \mathcal{P}_k} (s_p^0 +
    \sum\limits_{h\in\mathcal{K}}s_p^h\xi_h) \ge \bar d_k +\xi_k\hat d_k , ~k\in \mathcal{K},
    \xi\in \Xi,\label{con:example}\\
    & \sum\limits_{p\in \mathcal{P}_g} (s_p^0 +
    \sum\limits_{h\in\mathcal{K}}s_p^h\xi_h) -p_g^0 -
    \sum\limits_{h\in\mathcal{K}}p_g^h\xi_h  \leq 0,  g  \in G,\xi\in \Xi,
\end{align}
\begin{align}
    & \sum\limits_{k\in \mathcal{K}}\sum\limits_{p\in\mathcal{P}_k: l\in p}
    (s_p^0 + \sum\limits_{h\in\mathcal{K}}s_p^h\xi_h) \le \bar{y_l}^2R_lx_l,
    ~l \in E,\xi\in \Xi,\\
    &p_g^0 + \sum\limits_{h\in\mathcal{K}}p_g^h\xi_h \geq
    0, ~g\in G,~\xi\in \Xi, \\
    & s_p^0 + \sum\limits_{h\in\mathcal{K}}s_p^h \xi_h\geq 0, ~ p\in\mathcal{P}_k, ~ k  \in
    \mathcal{K}, \xi\in \Xi \\
    &x_l\in \{0,1\},~l\in E.
\end{align}
For each constraint, we can find its robust counterpart.
For example, the robust counterpart of constraint \eqref{con:example} writes
\begin{align}
 \sum\limits_{g\in G}  c_gp_g^0 + \sum\limits_{k \in\mathcal{K}}w^k\tau + w^0\kappa  \leq \gamma, \\
     w^k + w^0 \ge \sum\limits_{g\in G} p^k_g, k\in \mathcal{K}, \\
     w^k + w^0 \ge -\sum\limits_{g\in G} p^k_g, k\in \mathcal{K}, \\
     w^k \ge 0, w^0 \ge 0 , k\in \mathcal{K}.
\end{align}
See Appendix for the complete formulations w.r.t. uncertainty set $\Xi$ and
$\Xi'$. 

Also, note that the affinely adjustable formulation has a
two-fold solution: the solution to the line expansion decisions, and a set of
affine decision rules. The performance of adjustable decision rules will be
quantified by metrics on the optimal values of the affinely adjustable
robust-counterpart [PB2-AAR] and the exact robust adjustable counterpart [PB2-ARC].

%% file: dcopf-aar-results.tex
\section{Computational Illustrations}
\label{sec:results} 




{
\begin{table}[tbp]
\caption{Price of Robustness}
\label{tab:priceofrobustness}
\begin{center}
\begin{tabular}{l*{4}{c}}
\hline
\hline
 PB2-AAR$\Xi$ & Full protection & $\kappa$=3 & $\kappa$=2 & PB2 \\
\hline
Obj      & 1312 & 1288 & 1256 & 1160\\
Cons     & 500  & 500  & 500  & 26  \\
Vars     & 397  & 397  & 397  & 36  \\
Nodes    & 0    & 0    & 0    & 0   \\
Time(s.) & 0.19 & 0.15 & 0.19 & 0.02   \\
\hline
PB2-AAR$\Xi'$ & Full protection & $\alpha$=0.5 & $\alpha$=0 & PB2 \\
\hline
Obj      & 1332 & 1268 & 1217 & 1160\\
Cons     & 495  & 495  & 495  & 26\\
Vars     & 622  & 622  & 622  & 36  \\
Nodes    & 1    & 1    & 0    & 0   \\
Time(s.) & 0.28 & 0.21 & 0.17 & 0.02   \\
\hline                                   
\hline
\end{tabular}
\end{center}
\end{table}
}


\emph{Solvers}: 
We have formulated the models twice, once in AMPL, and once in using YALMIP.
In both cases, the paths were generated by Yen's K-shortest paths algorithm \cite{Kshortest} and the run-time of the path generation is not included
in the run-time of the solver.  
We have used IBM ILOG CPLEX for solving the AAR models 
and SeDuMi for solving the SDP relaxations of the AR models. 
In solving the semidefinite programs, we employ the sparsity-exploiting techniques
of Kim et al. \cite{KKMY}

\emph{Settings of parameters}: For the uncertainty set $\Xi$, 
Fix $\tau = 1$ to ensure
$\xi\in[-1,1]^{\abs{\mathcal{K}}}$. Particularly, for instance
Garver, set $\kappa\in \{ 2,3, \abs{\mathcal{K}}\}$. Note that when $\kappa =
\abs{\mathcal{K}}$, the fully protection is achieved, i.e., protection on set
$[-1,1]^{\abs{\mathcal{K}}}$. 
For the value-at-risk uncertainty set $\Xi'$, we set
the number of observations to be $N=10$ and $\alpha\in\{0,0.5, (N-1)/N\}$. And
uniformly generate $N$ values in $[\bar d_k-\hat d_k, \bar d_k - \hat d_k]$ for
each load bus $k\in \mathcal{K}$. Also note that
when $\alpha = (N-1)/N$, the value-at-risk uncertainty set is the convex hull
of all $N$ observations. 
To test the influences of the number of accessible paths to the solution, we
enumerate the value of $\abs{\mathcal{P}_{kg}}$ being $1,2,5$. 

\emph{Instances}: In our initial experiments, we have used a modification of the instance of Garver, 
due to Marecek et al. \cite{JonasPSCC}
As in Garver's original example we 
consider the connection of bus 6 to the existing system. 
In particular, we consider three double circuit 
lines 2-6, 3-6, and 4-6, and hence 8 
possible configurations. 
The existing lines are
complemented by an extra circuit along line 1-5 and 3-5 and 
the corresponding lines are replaced by their equivalent 
single circuit line. 
Line investment costs are uniformly at 100. 
See Figure \ref{fig:garver6instance} and Table~\ref{tab:1} for details.
In \cite{JonasPSCC}, it has been shown that only two configurations are AC feasible and piece-wise linearisations fail to provide the global optimum.

\emph{Results}:
In Table \ref{tab:priceofrobustness}, we summarise our initial results, including:
\begin{itemize}
    \item Obj: the transmission expansion cost and the power generation 
        cost.
    \item Cons: the number of linear constraints in affinely adjustable models.
    \item Var: the number of variables in affinely adjustable models. 
    \item Nodes: the number of branch and bound nodes by IBM ILOG CPLEX 12.5 with default parameters.  
    \item Time: CPU time in second.
    \end{itemize}
Albeit preminary, these results seem encouraging. 
We aim to provide a more comprehensive computational study in an extended version of the paper.

%% file: dcopf-aar-appendix.tex
\appendix
\vspace{-2em}
\section{An Appendix}\label{app:one}
The model associated with set $\Xi$ can be rewritten as 

\begin{align}
    \min& \hspace{0.5em}  \sum\limits_{l \in E }f_lx_l + \gamma
    \tag{PB2-AAR$\Xi$} \\
 \textrm{s.t.}& \hspace{0.5em} \sum\limits_{g\in G}  c_gp_g^0 +
    \sum\limits_{k \in\mathcal{K}}w^k\tau + w^0\kappa  \leq \gamma, \\
    & \hspace{0.5em} w^k + w^0 - \sum\limits_{g\in G} p^k_g \geq 0, k\in \mathcal{K}, \\
    & \hspace{0.5em} w^k + w^0 + \sum\limits_{g\in G} p^k_g \ge 0, k\in \mathcal{K}, \\
    & \hspace{0.5em} w^k \ge 0, w^0 \ge 0 , k\in \mathcal{K}, \\
    & \hspace{0.5em} \bar d_k - \sum\limits_{p\in \mathcal{P}_k} s^0_p
    + \sum\limits_{h\in \mathcal{K}} \tau u^h_k + \kappa u^0_k \leq 0, k\in \mathcal{K}, \\
    & \hspace{0.5em} u^h_k + u^0_k + \sum\limits_{p\in \mathcal{P}_k} s^h_p \ge
    0, h\neq k,  (h,k)\in
    \mathcal{K}^2,\\
    & \hspace{0.5em} u^h_k + u^0_k -  \sum\limits_{p\in \mathcal{P}_k} s^h_p
    \ge 0, h\neq k, (h,k)\in
    \mathcal{K}^2,\\
    & \hspace{0.5em} u^k_k + u^0_k - \hat d^k + \sum\limits_{p\in
\mathcal{P}_k} s^k_p \ge 0,  k\in
    \mathcal{K},\\
    &\hspace{0.5em}  u^k_k + u^0_k+\hat d^k - \sum\limits_{p\in \mathcal{P}_k}
    s^k_p \ge 0,  k\in
    \mathcal{K},\\
    &\hspace{0.5em}  u^h_k\geq 0, u^0_k\geq 0, (h,k)\in \mathcal{K},\\
    &\hspace{0.5em}  \sum\limits_{p \in \mathcal{P}_g} s^0_p - p^0_g +
    \sum\limits_{k \in \mathcal{K}} \tau v^k_g + \kappa v^0_g \leq 0, g \in G\\
    & \hspace{0.5em} v^h_g + v^0_g - \sum\limits_{p\in \mathcal{P}_g} s^h_p +
    p^h_g \ge 0, h \in \mathcal{K}, g\in G, \\
    & \hspace{0.5em} v^h_g + v^0_g  + \sum\limits_{p\in \mathcal{P}_g}
    s^h_p -  p^h_g  \ge 0, h \in \mathcal{K}, g\in G, \\
    & \hspace{0.5em} v^h_g \ge 0, v^0_g  \ge 0 , h \in \mathcal{K}, g\in G \\
 &\hspace{0.5em} \sum\limits_{k\in \mathcal{K}}\sum\limits_{p\in \mathcal{P}_k:
l\in p} s^0_p+\sum\limits_{k \in \mathcal{K}} \tau t^k_l + \kappa t^0_l  \leq  \bar y_l^2R_lx_l, l \in E, \\
    & \hspace{0.5em} t^h_l + t^0_l - \sum\limits_{k\in \mathcal{K}}
\sum\limits_{p\in \mathcal{P}_k:l \in p} s^h_p \ge 0 , h \in \mathcal{K}, l\in E, \\
    & \hspace{0.5em} t^h_l + t^0_l  + \sum\limits_{k\in \mathcal{K}}
\sum\limits_{p\in \mathcal{P}_k:l \in p} s^h_p \ge 0, h \in \mathcal{K}, l\in E, \end{align}

\begin{align}
    & \hspace{0.5em} t^h_l \ge 0, t^0_l  \ge 0 , h \in \mathcal{K}, l\in E \\
    &\hspace{0.5em}  \sum\limits_{h \in \mathcal{K}} \tau r^k_g + \kappa r^0_g  - p_g^0 \le 0, g\in G\\
    & \hspace{0.5em} r^h_g + r^0_g + p^h_g \ge 0 , h \in \mathcal{K}, g\in G, \\
    & \hspace{0.5em} r^h_g + r^0_g - p^h_g \ge 0, h \in \mathcal{K}, g\in G, \\
    & \hspace{0.5em} r^k_g \ge 0, r^0_g  \ge 0 , k \in \mathcal{K}, g\in G, \\
    &\hspace{0.5em}  \sum\limits_{h\in\mathcal{K}}\tau \mu^h_{kp}
+\kappa\mu^0_{kp} - s_p^0  \le 0 ,  p\in \mathcal{P}_k, k\in \mathcal{K} \\ &
\hspace{0.5em} s^h_p + \mu^h_{kp} + \mu^0_{kp} \ge 0, (h,k) \in
    \mathcal{K}^2, p\in \mathcal{P}_k, \\
    & \hspace{0.5em} \mu^h_{kp} + \mu^0_{kp} - s^h_p \ge 0 , (h,k) \in
    \mathcal{K}^2, p\in \mathcal{P}_k, \\
    & \hspace{0.5em} \mu^h_{kp} \ge 0, \mu^0_{kp}  \ge 0 , (h,k) \in
    \mathcal{K}^2,  p\in \mathcal{P}_k,\\
&\hspace{0.5em} x_l \in \{0,1\}, l\in E.
\end{align}

The robust counterpart associated with set $\Xi'$ is:
{
\begin{align}
    \min& \hspace{0.5em}  \sum\limits_{l \in E }f_lx_l + \gamma      \tag{PB2-AAR$\Xi'$} \\ 
    \textrm{s.t.}& \hspace{0.5em} \sum\limits_{g\in G}  c_gp_g^0 + w^0+
    \frac{1}{N(1-\alpha)}\sum\limits_{i=1}^N  w^i \leq \gamma,  \\
    &\hspace{0.5em} w^0 + w^i \geq \sum\limits_{g\in G} \sum\limits_{k\in
\mathcal{K}}p_g^k d_k^i  , i=1, \dots, N\\
    &\hspace{0.5em} w^i \ge 0 ,i=1, \dots, N  \\
& \hspace{0.5em}  u^0_k +   \frac{1}{N(1-\alpha)}\sum\limits_{i=1}^N u^i_k
-\sum\limits_{p\in \mathcal{P}_k} s^0_p  \le  0, k\in \mathcal{K}, \\
&\hspace{0.5em} u^0_k + u^i_k \geq d_k^i-\sum\limits_{p\in\mathcal{P}_k}\sum\limits_{h\in \mathcal{K}} s^h_p
d_h^i, i=1, \dots, N, k \in \mathcal{K}\\ &\hspace{0.5em} u^i_k \ge 0,
i=1,\dots,N,k \in \mathcal{K}\\
    &\hspace{0.5em}  \sum\limits_{p \in \mathcal{P}_g} s^0_p- p^0_g
 +v^0_g + \frac{1}{N(1-\alpha)}\sum\limits_{i=1}^N v^i_g    \leq 0,  g \in G,\\
 &\hspace{0.5em} v^0_g + v^i_g \geq \sum\limits_{p\in\mathcal{P}_g}\sum\limits_{h\in\mathcal{K}}s^h_pd_h^i -\sum\limits_{h\in \mathcal{K}}p_g^h d_h^i  , i=1, \dots, N, g \in G,\\
    &\hspace{0.5em} v^i_g \ge 0 ,i=1, \dots, N, g \in G
\end{align}

\begin{align}
    &\hspace{0.5em} \sum\limits_{k\in \mathcal{K}}\sum\limits_{p\in
    \mathcal{P}_k: l\in p} s^0_p+t^0_l + \frac{1}{N(1-\alpha)}\sum\limits_{i=1}^N t^i_l  \le \bar y_l^2R_lx_l,l \in E, \\
    &\hspace{0.5em} t^0_l + t^i_l \geq \sum\limits_{k\in
    \mathcal{K}}\sum\limits_{p\in \mathcal{P}_k: l\in p}\sum\limits_{h\in
    \mathcal{K}}  s^h_pd_h^i , i=1, \dots, N, l \in E,\\
    &\hspace{0.5em} t^i_l \ge 0 ,i=1, \dots, N, l \in E \\
    &\hspace{0.5em}  r^0_g + \frac{1}{N(1-\alpha)}\sum\limits_{i=1}^N
    r^i_g -p_g^0\le 0, g\in G\\
   &\hspace{0.5em} r^0_g + r^i_g  +\sum\limits_{h\in\mathcal{K}}p_g^hd_h^i \ge
    0, i=1, \dots, N,
   g \in G,\\
    &\hspace{0.5em} r^i_g \ge 0 ,i=1, \dots, N,  g \in G\\
    &\hspace{0.5em} \mu^0_{kp} + \frac{1}{N(1-\alpha)}\sum\limits_{i=1}^N
   \mu^i_{kp}  - s_p^0\le 0 , p\in \mathcal{P}_k, k\in \mathcal{K} \\
   &\hspace{0.5em} \mu^0_{kp} + \mu^i_{kp} +
   \sum\limits_{h\in\mathcal{K}}s_p^hd_h^i \ge 0, i=1, \dots, N, p\in \mathcal{P}_k,k \in \mathcal{K}\\
        &\hspace{0.5em} \mu^i_{kp} \ge 0 ,i=1, \dots, N,p\in \mathcal{P}_k, k \in \mathcal{K}\\
&\hspace{0.5em} x_l \in \{0,1\}, l\in E.
\end{align}
}
\vspace{-2em}